\documentclass[letterpaper, 10pt, conference]{ieeeconf}
\IEEEoverridecommandlockouts  
\usepackage{times} 
\usepackage{amsmath} 
\usepackage{amssymb}  
\usepackage{amsfonts}
\usepackage{mathrsfs}
\usepackage{graphicx}
\usepackage{xtab}
\usepackage{bbm}
\usepackage{bm}
\usepackage[euler]{textgreek}
\usepackage[utf8]{inputenc}
\usepackage[english]{babel}

\usepackage{xcolor}
\colorlet{GREY1}{white!20!black}
\colorlet{GREY2}{white!40!black}
\colorlet{GREY3}{white!60!black}

\usepackage[colorlinks=true,
            linkcolor=blue,
            urlcolor=blue,
            citecolor=blue,
            anchorcolor=blue,
            ]{hyperref}
\usepackage[capitalise]{cleveref}

\usepackage{tikz}
\usepackage{pgfplots}
\pgfplotsset{compat=1.10}
\usepgfplotslibrary{fillbetween}
\usetikzlibrary{pgfplots.groupplots}
\usetikzlibrary{matrix,shapes,arrows,positioning,chains}

\usetikzlibrary{calc}
\def\centerarc[#1](#2)(#3:#4:#5)
    { \draw[#1] ($(#2)+({#5*cos(#3)},{#5*sin(#3)})$) arc (#3:#4:#5); }
\tikzstyle{block} = [rectangle, draw, fill=white,
    text width=8em, text centered, rounded corners, minimum height=2em]
\tikzstyle{line} = [draw, -latex']

\makeatletter
\let\NAT@parse\undefined
\makeatother
\usepackage[numbers]{natbib}

\newcommand{\Rey}{\mbox{\textit{Re}}}
\newcommand{\Ca}{\mbox{\textit{Ca}}}

\newcommand{\norm}[1]{\left\lVert#1\right\rVert} 
\newcommand*{\tran}{^{\mkern-1.5mu\mathsf{T}}} 
\renewcommand{\matrix}[1]{
  \begin{pmatrix}
    #1
  \end{pmatrix}
} 

\makeatletter
\newcommand*{\diff}%
    {\@ifnextchar^{\DIfF}{\DIfF^{}}}
\def\DIfF^#1{%
    \mathop{\mathrm{\mathstrut d}}%
        \nolimits^{#1}\gobblespace}
\def\gobblespace{\futurelet\diffarg\opspace}
\def\opspace{%
    \let\DiffSpace\!%
    \ifx\diffarg(%
        \let\DiffSpace\relax
    \else
        \ifx\diffarg[%
            \let\DiffSpace\relax
        \else
            \ifx\diffarg\{%
                \let\DiffSpace\relax
            \fi\fi\fi\DiffSpace}

\title{\LARGE \bf
Stabilisation of falling liquid films with restricted observations*
}

\author{Oscar A. Holroyd\(^{1}\), Radu Cimpeanu\(^{1}\), and Susana N. Gomes\(^{1}\)
\thanks{*This work is supported by the UK Engineering and Physical Sciences Research Council (EPSRC) grant EP/S022848/1}%
\thanks{\(^{1}\)Warwick Mathematics Institute, University of Warwick, Coventry CV4 7AL, UK
        {\tt\small o.holroyd@warwick.ac.uk}}%
}

\begin{document}

\maketitle
\thispagestyle{empty}
\pagestyle{empty}

\begin{abstract}

We propose a method to stabilise a solution to equations describing the interface of thin liquid films falling under gravity with a finite number of actuators and restricted observations. As for many complex systems, full observation of the system state is challenging in physical settings, so methods able to take this into account are important. The Navier-Stokes equations modelling the flow are a complex, highly nonlinear set of PDEs, so standard control theoretical results are not applicable. Instead, we chain together a hierarchy of increasingly idealised approximations, developing a control strategy for the simplified model which is shown to be successfully applied to simulations of the full system.

\end{abstract}

\section{INTRODUCTION}
\label{sec:introduction}
  The stabilisation of liquid film interfaces is a prototypical example within a broader class of problems in which we seek to enable observation and control of complex physical systems, in this case the Navier-Stokes equations. Typical interfaces exhibit a broad range of parameter-dependent behaviours from travelling waves to chaos. Such flows have received much attention, both from experimental~\cite{nosoko1996characteristics} and analytical perspectives. Most pertinent to this work are a range of lower-dimensional long-wave models for the liquid-gas interface. These thin-film models, where the height of the film is assumed to be much smaller than its length, are comprehensively described in a review by \citeauthor{craster2009dynamics}~\cite{craster2009dynamics} and a book by \citeauthor{kalliadasis2011falling}~\cite{kalliadasis2011falling}. More recent work has included experimental validation~\cite{denner2018solitary} and additional asymptotic work to expand the range of validity of such models~\cite{usha2020evolution,mukhopadhyay2022modelling}.

  Typical applications of these films include coating flows~\cite{weinstein2004coating} --- requiring very flat interfaces --- and narrow channel micro-cooling~\cite{darabi2001electrohydrodynamic}, for which highly corrugated films are more desirable. Control capabilities targeted towards obtaining and maintaining different interfacial profiles are thus of key importance.

  Over the last ten years, progress has been made on producing feedback controllers for falling liquid films. \citeauthor{armaou2000feedback}~\cite{armaou2000feedback,christofides2000global} and \citeauthor{gomes2017controlling}~\cite{gomes2017controlling} applied control theoretical results from~\cite{zabczyk2020mathematical} to the Kuramoto-Sivashinsky (KS) equation, the simplest of the thin-film models, generating both computational and analytical results. Work by \citeauthor{thompson2016stabilising} showed that linear quadratic control (LQR) techniques are also effective for more complex models such as the Benney~\cite{benney1966long} and weighted-residual~\cite{ruyer2000improved} systems. Although \citeauthor{cimpeanu2021active}~\cite{cimpeanu2021active} demonstrated that simple proportional controls can successfully stabilise numerical simulations of the full Navier-Stokes equations, \citeauthor{holroyd2023linear}~\cite{holroyd2023linear} showed that a chain of increasingly accurate models can allow LQR techniques to be extended to the complete physical system, modelled by simulations of the Navier-Stokes equations.

  Previously we made the assumption of full observations~\cite{holroyd2023linear}, but here we introduce restrictions to a low-dimensional observation space, thus making an additional step towards physically realisable controls where interfacial observations are challenging~\cite{heining2013flow}. We describe a pair of control methodologies, previously applied to simpler systems in a study by \citeauthor{thompson2016stabilising}~\cite{thompson2016stabilising}, where observation of the interface is restricted to a finite number of discrete sites. We show that the same hierarchical control framework continues to permit successful stabilisation, albeit in a smaller subset of the parameter space.

\section{MATHEMATICAL MODELS}
\label{sec:governing_equations}
\begin{figure}[tpb]
    \centering
    \begin{tikzpicture}[font=\footnotesize]
      \draw [line width=0.25mm,-stealth] (0.25,0.2) -- +(0.10426666666, 0.658466666666);
      \draw [line width=0.25mm,-stealth] (0.25,0.2) -- +(0.658466666666,-0.10426666666);
      \draw (0.55,0.9) node{$y$};
      \draw (1.05,0.15) node{$x$};

      \draw [line width=0.25mm,-stealth] (6.35,1.4) -- +(0.0, -0.66666666666);
      \draw [line width=0.25mm,densely dashed] (5.69153333
    ,1.50426666) -- +(1.64616666,-0.2606666);
      \draw [line width=0.25mm] (5.69333,1.4) -- +(1.6666, 0.0);
      \centerarc[line width=0.25mm](6.35,1.4)(351:360:0.9)
      \draw (7.53333,1.3) node{$\theta$};
      \draw (6.55,0.85) node{$g$};

      \begin{axis}[
        width=1.1\columnwidth,
        axis lines = left,
        ymin = -5,
        ymax = 1.2,
        xmin = 0,
        xmax = 30,
        xticklabels={},
        yticklabels={},
        axis line style={draw=none},
        tick style={draw=none},
        unit vector ratio=1 1 1
        ]

        \addplot [
        name path = H,
        color = black,
        line width = 0.5mm,
        ]
        table[x=xh, y=h]
        {data_zhf.dat};

        \addplot [
        name path = B,
        color = black,
        line width = 0.5mm,
        dashed,
        ]
        table[x=xz, y=z]
        {data_zhf.dat};

        \addplot[black!20] fill between[of=H and B];

        \addplot [
        color = red,
        line width = 0.5mm,
        ]
        table[x=xf, y=f]
        {data_zhf.dat};

        \node[] at (axis cs: 12,0.6) {\textcolor{black}{$h(x,t)$}};
        \node[] at (axis cs: 12,-3) {\textcolor{black}{$f(x,t)$}};
      \end{axis}
    \end{tikzpicture}
    \caption{A thin-film with interface \(h\) under the influence of gravity \(g\) flows down a plane inclined an an angle \(\theta\). Discrete actuators injecting and removing fluid from the base (\(f\)) are used to stabilise the interface.}%
    \label{fig:diagram}
  \end{figure}
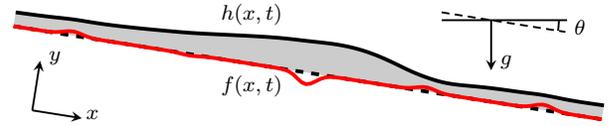

  As shown in \cref{fig:diagram}, we consider a thin film of fluid flowing down a plane inclined at an angle \(\theta\) relative to the horizontal. The fluid is described in a rotated coordinate system, so that \(y\) is normal to the base and \(x\) points downstream. On the wall at \(y=0\) we assume a no slip boundary condition with fluid injected and removed at a finite number of locations. At the upper boundary with the gas there is a free surface at \(y = h(x, t)\).
  For a film with mean height \(h_N\), the resulting gravitationally-driven uniform flow with flat interface at \(h = h_N\) (known as the Nusselt solution~\cite{nusselt1923warm}) has a parabolic velocity profile with surface velocity \(U_N\). For all but the most gentle regimes this solution is an unstable steady state, and so a control law is required to stabilise the flat interface.

  \subsection{Navier-Stokes Equations}
  \label{sub:navier_stokes_equations}
    For the incompressible fluids considered here, fluid flow can be modelled with the 2D Navier-Stokes equations. We assume that the density and viscosity ratios between the liquid film and the gas are sufficiently high that the effect of the gas flow on the liquid is negligible. This means that we can restrict our model to the liquid alone, which is parametrised by the Reynolds number \(\Rey\), capillary number \(\Ca\), domain length \(L\) and inclination angle \(\theta\). The liquid dynamics in the film is governed by the momentum equations
    \begin{align}
      \label{eqn:momentum1}\Rey (u_t + uu_x + vu_y) &= -p_x + 2 + u_{xx} + u_{yy}, \\
      \label{eqn:momentum2}\Rey (v_t + uv_x + vv_y) &= -p_y - 2\cot\theta + v_{xx} + v_{yy},
    \end{align}
    where \(u\) and \(v\) are the (dimensionless) horizontal and vertical velocity components and \(p\) the pressure, as well as the continuity equation
    \begin{equation}
      \label{eqn:continuity}u_x + v_y = 0.
    \end{equation}
    At the interface, \(y=h(x,t)\), we have the non-linear dynamic stress balance \cite{kalliadasis2011falling}
    \begin{align}
      \label{eqn:stress1}(v_x + u_y)(1-h_x^2) + 2h_x (v_y  -u_x) &= 0, \\
      \label{eqn:stress2}p - \frac{2(v_y + u_x h_x^2 - h_x(v_x + u_y))}{1+h_x^2} &= -\frac{1}{\Ca}\frac{h_{xx}}{(1+h_x^2)^{3/2}},
    \end{align}
    and no-slip and fluid injection/removal at the wall
    \begin{equation}
      \label{eqn:wall}u = 0, \qquad v = f(x,t).
    \end{equation}
    Finally the system is completed by the kinematic boundary condition
    \begin{equation}
      \label{eqn:kinematic}h_t = v-uh_x,
    \end{equation}
    and periodic boundaries in the \(x\)-direction.

    Defining the down-slope flux by integrating over the height of the film,
    \begin{equation}
      q(x,t) = \int_{0}^{h} u(x,y,t) \diff y,
    \end{equation}
    we combine \cref{eqn:continuity,eqn:wall,eqn:kinematic} to get the mass conservation equation
    \begin{equation}
      \label{eqn:conservation}h_t + q_x = f.
    \end{equation}

    The full Navier-Stokes system \cref{eqn:momentum1,eqn:momentum2,eqn:continuity,eqn:stress1,eqn:stress2,eqn:wall,eqn:kinematic} cannot in general be solved analytically, and so approximate solutions are found using a volume-of-fluids approach implemented by the free, open-source software Basilisk~\cite{basilisk,popinet2003}. All simulations are performed in a domain of width \(L=30\) inclined at an angle \(\theta=\pi/3\) with \(\Ca=0.05\).

  \subsection{Weighted-Residual Equations}
  \label{sub:weighted_residual_equations}
\begin{figure}[tpb]
      \centering
      \begin{tikzpicture}[font=\footnotesize]
        \begin{axis}[
            width=0.95\columnwidth,
            height=0.34\columnwidth,
            xlabel={\(x\)},
            ylabel={\(h\)},
            xmin=0, xmax=30,
            ymin=0, ymax=2,
            xtick={0, 5, 10, 15, 20, 25, 30},
            ytick={0, 1, 2},
        ]

        \node[right] at (axis cs: 0.2,1.5) {$t = 300$};
        \addplot [color = red, densely dashed,]
        table[x expr=\thisrowno{0}, y expr=\thisrowno{1}, skip first n=1]
        {data_wr_travelling_wave.dat};
        \addplot [color = black, solid,]
        table[x expr=\thisrowno{0}, y expr=\thisrowno{1}, skip first n=1]
        {data_ns_travelling_wave.dat};
        \end{axis}
      \end{tikzpicture}
      \caption{Fully evolved travelling waves for Navier solid-Stokes (in solid black) and weighted-residual (in dashed red) for \(\Rey = 10\). The weighted-residual equations provide a reasonable approximation to the Navier-Stokes interface, especially around the largest perturbation.}%
      \label{fig:interfaces}
    \end{figure}
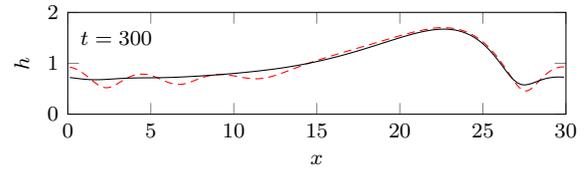

    \Cref{eqn:momentum1,eqn:momentum2,eqn:continuity,eqn:stress1,eqn:stress2,eqn:wall,eqn:kinematic} are complex and highly nonlinear, thus precluding the explicit construction of a control \(f\). In order to make further progress we use \citeauthor{ruyer2000improved}'s improved weighted-residual methodology~\cite{ruyer2000improved} which approximates \(u\) by a truncated sum of basis functions satisfying no-slip at the wall and zero tangential stress at the interface. Extended to include the forcing used here, the first-order truncation~\cite{thompson2016falling} gives the following expression for the flux
    \begin{equation}
    \begin{split}
      \label{eqn:wr}\frac{2\Rey}{5}h^2q_t + q ={}& \frac{h^3}{3}\left( 2 - 2 h_x \cot\theta + \frac{h_{xxx}}{\Ca} \right) \\
      &+ \Rey \left( \frac{18q^2h_x}{35} - \frac{34hqq_x}{35} + \frac{hqf}{5} \right).
    \end{split}
    \end{equation}
    Combined with \cref{eqn:conservation} this forms a 1D model for the evolution of the thin film. This reduced-dimensional system provides excellent agreement with the full Navier-Stokes equations up to \(\Rey \approx 5\)~\cite{ruyer2000improved}, significantly better than alternative models~\cite{salamon1994traveling}. \Cref{eqn:wr} captures many features of the Navier-Stokes systems including spontaneous back-flow~\cite{oron2009numerical}, and \cref{fig:interfaces} shows that some larger features are well-approximated at higher Reynolds numbers.

\section{CONTROL STRATEGIES}
\label{sec:control_strategies}
  To stabilise the Nusselt film, we focus on optimal feedback controls, in which the control is a function of observations of the interface. Denoting the system state (either \((u, v, p, \hat{h})\) or \((\hat{h}, q)\) --- where \(\hat{h} = h-1\) --- for Navier-Stokes or weighted-residual systems respectively) by \(\xi \in X\), where \(X\) is an appropriate Hilbert space (e.g. \(h \in L^2(0, T; H^2(0, L))\), \(q \in L^2(0, T; H^1(0, L))\) with periodic boundaries), we can write the feedback control problem in an abstract sense:
  \begin{equation}
    \dot{\xi} = \mathcal{A}\xi + \mathcal{B}\eta, \quad \zeta = \mathcal{C}\xi, \quad \eta = \mathcal{K}\zeta,
  \end{equation}
  where \(\eta \in \mathbb{R}^M\) is the control, \(\zeta \in \mathbb{R}^P\) is some system observation and \(\mathcal{A}\), \(\mathcal{B}\), \(\mathcal{C}\), and \(\mathcal{K}\) are operators acting between relevant spaces: the uncontrolled PDE operator, the control actuator mechanism, the observation operator, and the gain operator to be computed. This can be written in closed-loop form
  \begin{equation}
    \label{eqn:closed_loop}\dot{\xi} = (\mathcal{A} + \mathcal{B}\mathcal{K}\mathcal{C})\xi,
  \end{equation}
  which makes it clear that \(\mathcal{K}\) must be chosen so that the operator \(\mathcal{A} + \mathcal{B}\mathcal{K}\mathcal{C}\) is stable. Finally, to quantify the performance of a given gain operator \(\mathcal{K}\), we define the cost
  \begin{equation}
    \label{eqn:cost}\kappa = \int_{0}^{\infty} \int_{0}^{L}  \beta \hat{h}(x)^2 + (1-\beta) f^2 \diff x \diff t.
  \end{equation}

  Thus far we have not restricted the form of our control function \(f\). However, arbitrary injection and removal of fluid along a continuum of locations at the base is physically unfeasible, and so we restrict our controls to a finite number of injection sites \(\{x_i\}_{i=1}^M\) (and so the control space is \(\mathbb{R}^M\)), i.e.
  \begin{equation}
    f(x, t) = \sum_{i=1}^{M} \eta_i(t) d(x-x_i),
  \end{equation}
  where \(\eta_i(t)\) are the individual, time-dependent, control amplitudes, and \(d\) is a smooth, periodic approximation to the Dirac-delta distribution
  \begin{equation}
    d(x) = \alpha \exp\left[ \frac{\cos(2\pi x/L)-1}{\omega^2} \right],
  \end{equation}
  where \(\omega\) controls the width of the function, and \(\alpha\) is chosen so that \(\int_{0}^{L} d(x) \diff x = 1\).

  Despite the simplifying assumptions made to convert the Navier-Stokes system to the weighted-residual system, the model retains high-order derivatives and significant nonlinearities. For such a PDE, it is not yet possible to explicitly compute the optimal \(\mathcal{K}\) because of these nonlinearities. In order to solve the problem, we must both linearise and spatially discretise the problem so we are left with a system of \(N\) linear ODEs (where \(N\) is the number of nodes in the discretisation):
  \begin{equation}
    \label{eqn:closed_loop_odes}\dot{\xi} = (A + BKC)\xi.
  \end{equation}
  The cost \(\kappa\) is similarly discretised to give
  \begin{equation}
    \label{eqn:discrete_cost}c = \int_{0}^{\infty} \hat{h}\tran U \hat{h} + \eta\tran V \eta \diff t,
  \end{equation}
  where \(U = \frac{\beta L}{N} I\in\mathbb{R}^{N\times N}\) and \(V = (1-\beta) I\in\mathbb{R}^{M\times M}\) are matrices whose entries are chosen as to form the discrete analogue of the continuous cost \cref{eqn:cost}.

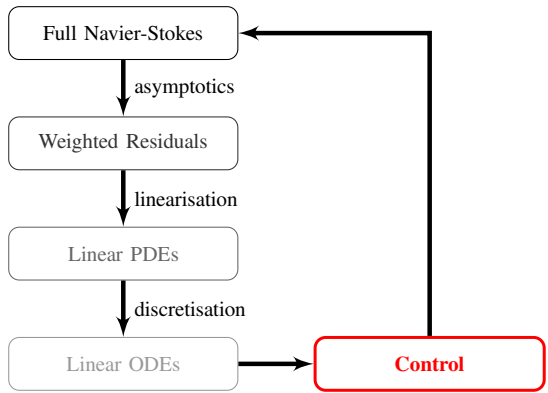
\begin{figure}[tpb]
    \centering
    \vspace{0.5em}
    \begin{tikzpicture}[font=\footnotesize, node distance = 0.75cm and 1cm, auto]
      \node [block] (one) {Full Navier-Stokes};
      \node [GREY1, block, below=of one] (two) {Weighted Residuals};
      \node [GREY2, block, below=of two] (three) {Linear PDEs};
      \node [GREY3, block, below=of three] (four) {Linear ODEs};
      \node [red, block, right=of four, line width=0.4mm] (five) {\textbf{Control}};

      \path [line,line width=0.6mm] (one) -- (two) node [midway, right] {asymptotics};
      \path [line,line width=0.6mm] (two) -- (three) node [midway, right] {linearisation};
      \path [line,line width=0.6mm] (three) -- (four) node [midway, right] {discretisation};
      \path [line,line width=0.6mm] (four) -- (five);
      \path [line,line width=0.6mm] (five) |- (one);
    \end{tikzpicture}
    \caption{Illustration of the chain of simplifications that permit the design of stabilising controls for the full Navier-Stokes simulation of a falling liquid film: by making a chain of increasingly abstracting simplifying assumptions we eventually reach a system where we can apply established control theoretic results. By projecting the full system to the space in which the controls are designed and then reprojecting the controls back to the full space we can attempt to stabilise solutions to the original problem.}%
    \label{fig:hierarchy}
  \end{figure}

  The stabilisation of systems of linear ODEs is a well-studied problem~\cite{zabczyk2020mathematical}, and so progress can now be made in computing the gain matrix \(K\). \citeauthor{al2018linearized}~\cite{al2018linearized} proved that such discrete approximations converge to a stabilising control for the KS equation. This ends the hierarchy of simplifications required to stabilise the Nusselt solution of the Navier-Stokes equations (summarised in \cref{fig:hierarchy}). We can test the methodology on direct numerical simulations of the Navier-Stokes film: we choose a discretisation space and compute the matrix \(K\) for the approximating linear system \cref{eqn:closed_loop_odes}, and, by applying it to observations projected into the discretisation space, we can use it as a feedback operator to compute the actuator strengths \(\eta_i\), which can be applied to observations of the Navier-Stokes simulation.

  \subsection{Full Observations (State Feedback LQR)}
  \label{sub:full_observations}
    In the case of full observations, when \(\mathcal{C}\) is the identity and \(C = I \in \mathbb{R}^{N\times N}\) finding the optimal \(K\) becomes a state-feedback linear quadratic regulator (LQR) problem, for which there is a wealth of literature~\cite{johnson1970design,syrmos1997static}. \citeauthor{holroyd2023linear}~\cite{holroyd2023linear} demonstrated that, in this case, such a hierarchical scheme provides a sufficiently robust approximation to successfully control the film. Linear stability analysis and numerical experiments suggest that a sufficient condition for stabilisation was that \(M\) should be at least as large as the number of unstable modes of the linearised system with equidistantly placed actuators. We note that it is not a necessary condition however --- numerous cases where fewer actuators are able to stabilise the uniform film were observed (see Fig. 8 in~\cite{holroyd2023linear}).

  \subsection{Restricted Observations (Output Feedback LQR)}
  \label{sub:restricted_observations}
    The control scheme devised by \citeauthor{holroyd2023linear}~\cite{holroyd2023linear} has access to full information about the interface. However, although this is trivial in a simulation, in reality it is almost impossible to make real-time measurements of the entire interface~\cite{heining2013flow}, rendering this method impractical for experiments. It is more realistic to restrict the observations of the interfacial height \(h\) to a finite-dimensional space corresponding to discrete observations from \(P\) arbitrarily placed observers, rather than the full state. This corresponds to \(C \in \mathbb{R}^{P\times N}\).

    In the periodic case here, when we are able to observe the entire interface we can use the translational invariance of the problem to argue that the optimal placement of the actuators is evenly spaced throughout the domain~\cite{tomlin2019point}. The optimal placement of observers and actuators for arbitrary pairs of \(M\) and \(P\) is a complicated problem~\cite{tomlin2019point} and is left for future work --- here we evenly space both actuators and observers around the centre of the domain.

    When \(P < N\) and \(C \neq I\), choosing the optimal \(K\) to minimise \cref{eqn:discrete_cost} is much more challenging than the basic LQR problem in \cref{sub:full_observations}. Finding the matrix necessitates solving the static output feedback problem
    \begin{align}
      \label{eqn:cond1}0 &= A_c\tran Q + Q A_c + U + C\tran K\tran V K C, \\
      \label{eqn:cond2}0 &= A_c S + S A_c\tran + I, \\
      \label{eqn:cond3}0 &= VKCSC\tran + B\tran QSC\tran
    \end{align}
    for \(K\) (as well as \(Q\) and \(S\)), where \(A_c = A+BKC\). In the case when \(C = I\) \cref{eqn:cond1,eqn:cond2,eqn:cond3} collapse to the continuous algebraic Riccati equation which can be solved directly~\cite{levine1970determination}, but the general problem must be solved using expensive iterative methods~\cite{ilka2023novel}.

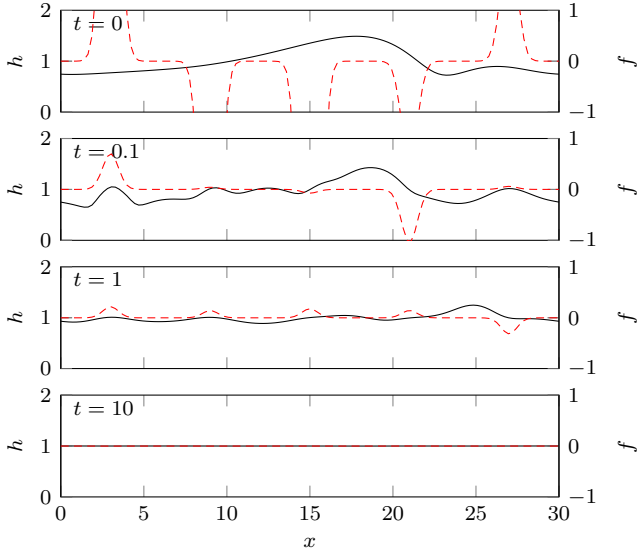
\begin{figure}[tpb]
      \centering
      \begin{tikzpicture}[font=\footnotesize]
        \begin{groupplot}[
          group style={
              group name=my plots,
              group size=1 by 4,
              xlabels at=edge bottom,
              xticklabels at=edge bottom,
              ylabels at=edge left,
              yticklabels at=edge left,
              vertical sep=10pt,
              horizontal sep=10pt
          },
          width=0.95\columnwidth,
          height=0.34\columnwidth,
          xlabel={\(x\)},
          ylabel={\(h\)},
          xmin=0, xmax=30,
          ymin=0, ymax=2,
          xtick={0, 5, 10, 15, 20, 25, 30},
          ytick={0, 1, 2},
          ]
          \nextgroupplot
          \def \num {100}
          \node[right] at (axis cs: 0.2,1.75) {$t = 0$};
          \addplot [color = black, solid,]
          table[x expr=\thisrowno{0}, y expr=\thisrowno{1}, skip first n=1]
          {data_sof-out_ns-1-0000000\num.dat};

          \nextgroupplot
          \def \num {101}
          \node[right] at (axis cs: 0.2,1.75) {$t = 0.1$};
          \addplot [color = black, solid,]
          table[x expr=\thisrowno{0}, y expr=\thisrowno{1}, skip first n=1]
          {data_sof-out_ns-1-0000000\num.dat};

          \nextgroupplot
          \def \num {110}
          \node[right] at (axis cs: 0.2,1.75) {$t = 1$};
          \addplot [color = black, solid,]
          table[x expr=\thisrowno{0}, y expr=\thisrowno{1}, skip first n=1]
          {data_sof-out_ns-1-0000000\num.dat};

          \nextgroupplot
          \def \num {200}
          \node[right] at (axis cs: 0.2,1.75) {$t = 10$};
          \addplot [color = black, solid,]
          table[x expr=\thisrowno{0}, y expr=\thisrowno{1}, skip first n=1]
          {data_sof-out_ns-1-0000000\num.dat};
        \end{groupplot}
        \begin{groupplot}[
          group style={
              group name=my plots,
              group size=1 by 4,
              xlabels at=edge bottom,
              xticklabels at=edge bottom,
              ylabels at=edge right,
              yticklabels at=edge right,
              vertical sep=10pt,
              horizontal sep=10pt
          },
          width=0.95\columnwidth,
          height=0.34\columnwidth,
          xlabel={\(x\)},
          ylabel={\(f\)},
          xmin=0, xmax=30,
          ymin=-1, ymax=1,
          ytick={-1, 0, 1},
          axis x line=none,
          ylabel near ticks,
          y label style={at={(axis description cs:1.1,0.5)}},
          ]
          \nextgroupplot
          \def \num {100}
          \node[right] at (axis cs: 0.2,1.75) {$t = 0$};
          \addplot [color = red, densely dashed,]
          table[x expr=\thisrowno{0}, y expr=\thisrowno{2}, skip first n=1]
          {data_sof-out_ns-1-0000000\num.dat};

          \nextgroupplot
          \def \num {101}
          \node[right] at (axis cs: 0.2,1.75) {$t = 0.5$};
          \addplot [color = red, densely dashed,]
          table[x expr=\thisrowno{0}, y expr=\thisrowno{2}, skip first n=1]
          {data_sof-out_ns-1-0000000\num.dat};

          \nextgroupplot
          \def \num {110}
          \node[right] at (axis cs: 0.2,1.75) {$t = 5$};
          \addplot [color = red, densely dashed,]
          table[x expr=\thisrowno{0}, y expr=\thisrowno{2}, skip first n=1]
          {data_sof-out_ns-1-0000000\num.dat};

          \nextgroupplot
          \def \num {200}
          \node[right] at (axis cs: 0.2,1.75) {$t = 50$};
          \addplot [color = red, densely dashed,]
          table[x expr=\thisrowno{0}, y expr=\thisrowno{2}, skip first n=1]
          {data_sof-out_ns-1-0000000\num.dat};
        \end{groupplot}
      \end{tikzpicture}
      \caption{Interface \(h\) in solid black, controls \(f\) in dashed red for increasing times after the actuators are switched on at \(t=0\) with \(\Rey = 11.29\) and  \(M=P=5\). Once the controls are switched on at the uniform film \(h=1\) is stabilised.}%
      \label{fig:sof_damping_film}
    \end{figure}

\begin{figure}[tpb]
      \centering
      \begin{tikzpicture}[font=\footnotesize]
        \begin{axis}[
          width=0.95\columnwidth,
          height=0.6\columnwidth,
          axis lines = box,
          xlabel = {\(t\)},
          ylabel = {\(\norm{h-1}\)},
          xmin = -50,
          xmax = 100,
          ymode=log,
          ]

          \draw [color = black, dotted,]
          ({axis cs:0,0}|-{rel axis cs:0,0}) -- ({axis cs:0,0}|-{rel axis cs:0,1});

          \addplot [color = black, solid,]
          table[x expr=\thisrowno{0}, y expr=\thisrowno{1}]
          {data_sof-out_ns-0.dat};

          \addplot [color = black, densely dashed, domain=0:100]
          {exp(-0.054485477865575*x)};
        \end{axis}
      \end{tikzpicture}
      \caption{Deviation of the interface from \(h=1\) for the regime in \cref{fig:sof_damping_film} (solid line). The perturbations are exponentially damped with rate exceeding that predicted by the linearised approximation to the dynamics (dashed line).}%
      \label{fig:sof_damping_norm}
    \end{figure}
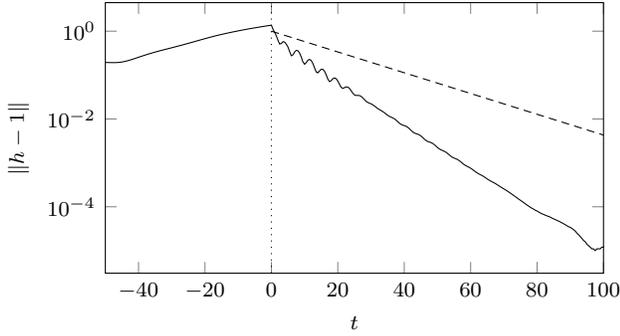

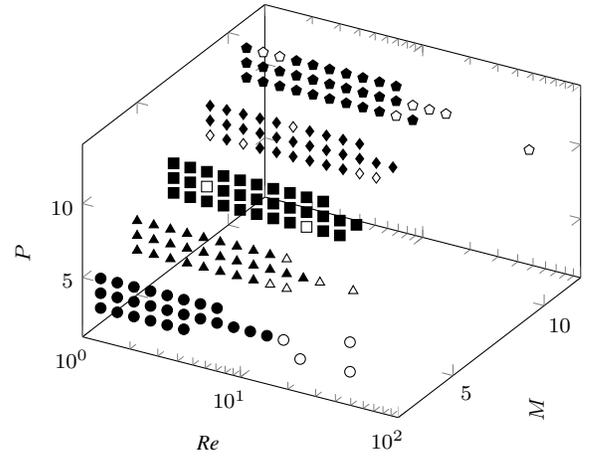
\begin{figure}[tpb]
      \centering
      \vspace{0.5em}
      \begin{tikzpicture}[font=\footnotesize]
        \begin{axis}[
          width=0.95\columnwidth,
          axis lines = box,
          xlabel = {\(\Rey\)},
          ylabel = {\(M\)},
          zlabel = {\(P\)},
          xmin = 1,
          xmax = 100,
          ymin = 2,
          ymax = 12,
          zmin = 1,
          zmax = 14,
          xmode=log,
          view={30}{40},
          ]

          \addplot3 [only marks, fill=black, mark=*,
            y filter/.expression={ y-3 == 0 ? \pgfmathresult : NaN)}]
          table[x=R, y=M, z=P]
          {data_sof-success.dat};

          \addplot3 [only marks, fill=white, mark=*,
            y filter/.expression={ y-3 == 0 ? \pgfmathresult : NaN)}]
          table[x=R, y=M, z=P]
          {data_sof-failure.dat};

          \addplot3 [only marks, fill=black, mark=triangle*,
            y filter/.expression={ y-5 == 0 ? \pgfmathresult : NaN)}]
          table[x=R, y=M, z=P]
          {data_sof-success.dat};

          \addplot3 [only marks, fill=white, mark=triangle*,
            y filter/.expression={ y-5 == 0 ? \pgfmathresult : NaN)}]
          table[x=R, y=M, z=P]
          {data_sof-failure.dat};

          \addplot3 [only marks, fill=black, mark=square*,
            y filter/.expression={ y-7 == 0 ? \pgfmathresult : NaN)}]
          table[x=R, y=M, z=P]
          {data_sof-success.dat};

          \addplot3 [only marks, fill=white, mark=square*,
            y filter/.expression={ y-7 == 0 ? \pgfmathresult : NaN)}]
          table[x=R, y=M, z=P]
          {data_sof-failure.dat};

          \addplot3 [only marks, fill=black, mark=diamond*,
            y filter/.expression={ y-9 == 0 ? \pgfmathresult : NaN)}]
          table[x=R, y=M, z=P]
          {data_sof-success.dat};

          \addplot3 [only marks, fill=white, mark=diamond*,
            y filter/.expression={ y-9 == 0 ? \pgfmathresult : NaN)}]
          table[x=R, y=M, z=P]
          {data_sof-failure.dat};

          \addplot3 [only marks, fill=black, mark=pentagon*,
            y filter/.expression={ y-11 == 0 ? \pgfmathresult : NaN)}]
          table[x=R, y=M, z=P]
          {data_sof-success.dat};

          \addplot3 [only marks, fill=white, mark=pentagon*,
            y filter/.expression={ y-11 == 0 ? \pgfmathresult : NaN)}]
          table[x=R, y=M, z=P]
          {data_sof-failure.dat};
        \end{axis}
      \end{tikzpicture}
      \caption{Successful (in black) and unsuccessful (in white) results for the output-feedback LQR method applied to simulations of the Navier-Stokes equations at various \(\Rey\), number of actuators \(M\), and numbers of observers \(P\). The groups of differently shaped points correspond to the values of \(M\): 3, 5, 7, 9, and 11. Gaps occur when the iterative process to find \(K\) failed.}%
      \label{fig:sof_damping_scatter}
    \end{figure}

    Even with restricted observations we are able to stabilise the uniform film solution in some cases, as shown in \cref{fig:sof_damping_film,fig:sof_damping_norm}.
    However, with observations restricted to \(\mathbb{R}^P\), the output feedback LQR control is no longer able to stabilise the film for large values of \(\Rey\) when the more significant inertial effects mean that the perturbations become harder to control. For large \(\Rey\) the difference between Navier-Stokes and weighted-residual equations is larger, and so sufficiently large waves can develop in the gaps between observers so as to make stabilisation impossible. In addition to this, in many cases the iterative algorithm in~\cite{ilka2023novel} either fails to start or fails to converge, likely due to the matrix problem being ill-posed thanks to insufficient observation of the system. \Cref{fig:sof_damping_scatter} shows how damping is possible for lower \(\Rey\) but for larger values either the control fails or a stabilising \(K\) cannot be found for \cref{eqn:closed_loop_odes}. Nevertheless, it is worth remarking that, as for static-output LQR~\cite{holroyd2023linear}, the framework functions beyond the Reynolds numbers where the predictive capacity of the weighted-residual model breaks down (albeit not as far beyond this limit as for full state observations).

  \subsection{Luenberger Observer}
  \label{sub:luenberger_observer}
    Although the static output feedback law can stabilise the uniform film as long as \(M\) and \(P\) are sufficiently large, it does not make maximum use of the available information. In the fully observed case, the observation at a given time \(t\) fully describes the state of the system, but there are a number of possible solutions that could be constructed to fit the \(P\) observations. We can improve upon this by allowing the use of past data to further our ability to construct the interface at the current time. For example, since the interface is generally convected downstream by the flux, it is reasonable to suggest that, shortly after an observation at \((t, x_i)\), the interfacial height at \((t + \delta t, x_i + \delta x)\) is likely to be similar. We can formalise this ability to make predictions about the current interface from past data and knowledge of the underlying dynamics by introducing a dynamic estimator \(z\), which is continually updated with the observation data but also allowed to evolve in time under its own separate dynamics. The method is described by \citeauthor{luenberger1971introduction}~\cite{luenberger1971introduction} and has been successfully applied directly to the KS equation~\cite{armaou2000feedback}, and other models including the Benney equation~\cite{thompson2016stabilising}. Here we use the same model, where we update our predictions with both the linearised system dynamics and the observations, but now as part of the hierarchical control framework.

    While our controls are part of a system in real space, there is no reason for the estimator to be. After transforming to Fourier space (so \(\tilde{\xi} = \mathcal{F}\hat{\xi}\) where \(\mathcal{F}\) is the Fourier transform), we reorder the wavenumbers to separate stable and unstable modes:
    \begin{equation}
      \dot{\tilde{\xi}} = \tilde{A}\tilde{\xi} + \tilde{B}\tilde{K}\tilde{\xi} =
      \matrix{\tilde{A}_u & 0 \\ 0 & \tilde{A}_s}\tilde{\xi} + \matrix{\tilde{B}_u \\ \tilde{B}_s}\tilde{K}\tilde{\xi}.
    \end{equation}
    Separating the unstable modes more explicitly, we write
    \begin{equation}
    \begin{split}
      \matrix{\dot{\tilde{\xi}}_u \\ \dot{\tilde{\xi}}_s} &= \matrix{\tilde{A}_u & 0 \\ 0 & \tilde{A}_s}\matrix{\tilde{\xi}_u \\ \tilde{\xi}_s} + \matrix{\tilde{B}_u \\ \tilde{B}_s}\tilde{K}\matrix{\tilde{\xi}_u \\ \tilde{\xi}_s} \\
      &= \matrix{\tilde{A}_u + \tilde{B}_u \tilde{K} & 0 \\ \tilde{B}_s \tilde{K} & \tilde{A}_s}\matrix{\tilde{\xi}_u \\ \tilde{\xi}_s},
    \end{split}
    \end{equation}
    and we see that, since the matrix on the right hand side is block lower triangular, the controls leave the eigenvalues of the stable modes unchanged, and so they remain stable. We thus reduce the control problem to
    \begin{equation}
      \dot{\tilde{\xi}}_u = \tilde{A}_u \tilde{\xi}_u + \tilde{B}_u \tilde{K} \tilde{\xi}_u,
    \end{equation}
    for which we select an optimal \(\tilde{K}\) using the LQR algorithm. Since we are unable to directly measure \(\tilde{\xi}_u\), we replace it with the estimator \(z\):
    \begin{equation}
      \label{eqn:DOF}\dot{\tilde{\xi}}_u = \tilde{A}_u \tilde{\xi}_u + \tilde{B}_u \tilde{K} z.
    \end{equation}
    We note that there is an implicit assumption here that \(M\) --- the dimension of \(z\) --- must equal (or exceed) the number of unstable modes for the system to be controllable, and the same is true for \(P\) and observability.

    The dynamics of the estimator are governed by the same linearised dynamics as \(\tilde{\xi}_u\), with an additional forcing term to push the system to match the observations. This gives us a set of \(M\) ODEs,
    \begin{equation}
      \label{eqn:estimator}\dot{z} = (\tilde{A}_u + \tilde{B}_u \tilde{K}) z + L(C \xi - C \mathcal{F}_u^{-1} z),
    \end{equation}
    where \(\mathcal{F}_u^{-1}\) is the inverse Fourier transform mapping the \(M\) estimated modes to \(N\) real values.

    If we define the error \(e = \tilde{\xi}_u-z\), by substituting in \cref{eqn:DOF,eqn:estimator} and noting that \(h = \mathcal{F}_u^{-1}\tilde{\xi}_u + \mathcal{F}_s^{-1}\tilde{\xi}_s\) we obtain
    \begin{equation}
      \dot{e} = (\tilde{A}_u - LC\mathcal{F}_u^{-1}) e - LC\mathcal{F}_s^{-1}\tilde{\xi}_s.
    \end{equation}
    The stable terms will decay on their own and so all that remains is to choose \(L\) so that \(\tilde{A}_u - LC\mathcal{F}_u^{-1}\) is stable (i.e.\ all its eigenvalues have negative real part). We do this by noting that eigenvalues are invariant under transposition and use the LQR algorithm on \((\tilde{A}_u - LC\mathcal{F}_u^{-1})\tran = \tilde{A}_u\tran - (C\mathcal{F}_u^{-1})\tran L\tran\).

    As we can see in \cref{fig:dof_damping_film,fig:dof_damping_norm}, accurate information about the precise shape of the interface is not always necessary; it is sufficient to be able to observe the most unstable modes of the system.
    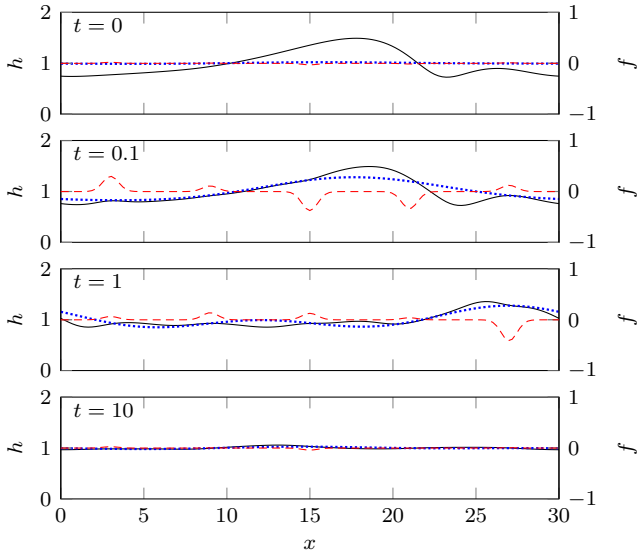
\begin{figure}[tpb]
      \centering
      \begin{tikzpicture}[font=\footnotesize]
        \begin{groupplot}[
          group style={
              group name=my plots,
              group size=1 by 4,
              xlabels at=edge bottom,
              xticklabels at=edge bottom,
              ylabels at=edge left,
              yticklabels at=edge left,
              vertical sep=10pt,
              horizontal sep=10pt
          },
          width=0.95\columnwidth,
          height=0.34\columnwidth,
          xlabel={\(x\)},
          ylabel={\(h\)},
          xmin=0, xmax=30,
          ymin=0, ymax=2,
          xtick={0, 5, 10, 15, 20, 25, 30},
          ytick={0, 1, 2},
          ]
          \nextgroupplot
          \def \num {100}
          \node[right] at (axis cs: 0.2,1.75) {$t = 0$};
          \addplot [color = black, solid,]
          table[x expr=\thisrowno{0}, y expr=\thisrowno{1}, skip first n=1]
          {data_dof-out_ns-1-0000000\num.dat};
          \addplot [color = blue, densely dotted, line width=0.3mm]
          table[x expr=\thisrowno{0}, y expr=1+\thisrowno{3}, skip first n=1]
          {data_dof-out_ns-1-0000000\num.dat};

          \nextgroupplot
          \def \num {101}
          \node[right] at (axis cs: 0.2,1.75) {$t = 0.1$};
          \addplot [color = black, solid,]
          table[x expr=\thisrowno{0}, y expr=\thisrowno{1}, skip first n=1]
          {data_dof-out_ns-1-0000000\num.dat};
          \addplot [color = blue, densely dotted, line width=0.3mm]
          table[x expr=\thisrowno{0}, y expr=1+\thisrowno{3}, skip first n=1]
          {data_dof-out_ns-1-0000000\num.dat};

          \nextgroupplot
          \def \num {110}
          \node[right] at (axis cs: 0.2,1.75) {$t = 1$};
          \addplot [color = black, solid,]
          table[x expr=\thisrowno{0}, y expr=\thisrowno{1}, skip first n=1]
          {data_dof-out_ns-1-0000000\num.dat};
          \addplot [color = blue, densely dotted, line width=0.3mm]
          table[x expr=\thisrowno{0}, y expr=1+\thisrowno{3}, skip first n=1]
          {data_dof-out_ns-1-0000000\num.dat};

          \nextgroupplot
          \def \num {200}
          \node[right] at (axis cs: 0.2,1.75) {$t = 10$};
          \addplot [color = black, solid,]
          table[x expr=\thisrowno{0}, y expr=\thisrowno{1}, skip first n=1]
          {data_dof-out_ns-1-0000000\num.dat};
          \addplot [color = blue, densely dotted, line width=0.3mm]
          table[x expr=\thisrowno{0}, y expr=1+\thisrowno{3}, skip first n=1]
          {data_dof-out_ns-1-0000000\num.dat};
        \end{groupplot}
        \begin{groupplot}[
          group style={
              group name=my plots,
              group size=1 by 4,
              xlabels at=edge bottom,
              xticklabels at=edge bottom,
              ylabels at=edge right,
              yticklabels at=edge right,
              vertical sep=10pt,
              horizontal sep=10pt
          },
          width=0.95\columnwidth,
          height=0.34\columnwidth,
          xlabel={\(x\)},
          ylabel={\(f\)},
          xmin=0, xmax=30,
          ymin=-1, ymax=1,
          ytick={-1, 0, 1},
          axis x line=none,
          ylabel near ticks,
          y label style={at={(axis description cs:1.1,0.5)}},
          ]
          \nextgroupplot
          \def \num {100}
          \node[right] at (axis cs: 0.2,1.75) {$t = 0$};
          \addplot [color = red, densely dashed,]
          table[x expr=\thisrowno{0}, y expr=\thisrowno{2}, skip first n=1]
          {data_dof-out_ns-1-0000000\num.dat};

          \nextgroupplot
          \def \num {101}
          \node[right] at (axis cs: 0.2,1.75) {$t = 0.5$};
          \addplot [color = red, densely dashed,]
          table[x expr=\thisrowno{0}, y expr=\thisrowno{2}, skip first n=1]
          {data_dof-out_ns-1-0000000\num.dat};

          \nextgroupplot
          \def \num {110}
          \node[right] at (axis cs: 0.2,1.75) {$t = 5$};
          \addplot [color = red, densely dashed,]
          table[x expr=\thisrowno{0}, y expr=\thisrowno{2}, skip first n=1]
          {data_dof-out_ns-1-0000000\num.dat};

          \nextgroupplot
          \def \num {200}
          \node[right] at (axis cs: 0.2,1.75) {$t = 50$};
          \addplot [color = red, densely dashed,]
          table[x expr=\thisrowno{0}, y expr=\thisrowno{2}, skip first n=1]
          {data_dof-out_ns-1-0000000\num.dat};
        \end{groupplot}
      \end{tikzpicture}
      \caption{The Luenberger estimator captures the most important modes of the interface, enabling the controls to stabilise the uniform film solution. Interface \(h\) in black, estimator \(z\) in dotted blue, controls \(f\) in dashed red. \(\Rey = 11.29\), \(M=P=5\).}%
      \label{fig:dof_damping_film}
    \end{figure}
    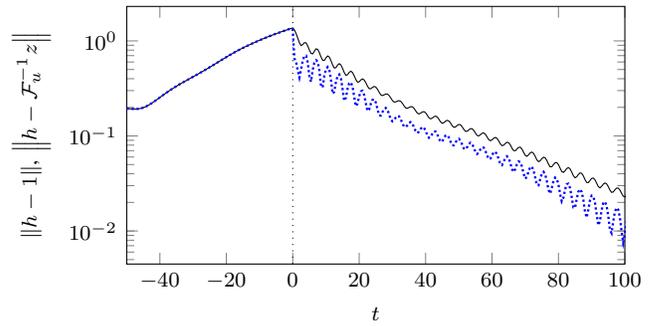
\begin{figure}[tpb]
      \centering
      \vspace{0.5em}
      \begin{tikzpicture}[font=\footnotesize]
        \begin{axis}[
          width=0.95\columnwidth,
          height=0.58\columnwidth,
          axis lines = box,
          xlabel = {\(t\)},
          ylabel = {\(\norm{h-1}\), \(\norm{h-\mathcal{F}_u^{-1}z}\)},
          xmin = -50,
          xmax = 100,
          ymode=log,
          ]

          \draw [color = black, dotted,]
          ({axis cs:0,0}|-{rel axis cs:0,0}) -- ({axis cs:0,0}|-{rel axis cs:0,1});

          \addplot [color = black, solid,]
          table[x expr=\thisrowno{0}, y expr=\thisrowno{1}]
          {data_dof-out_ns-0.dat};

          \addplot [color = blue, densely dotted, line width=0.3mm]
          table[x expr=\thisrowno{0}, y expr=\thisrowno{2}]
          {data_dof-out_ns-0.dat};
        \end{axis}
      \end{tikzpicture}
      \caption{Once the controls are activated, both the interfacial perturbations and estimator error (solid black and dotted blue respectively) decay exponentially. Parameters as for \cref{fig:dof_damping_film}.}%
      \label{fig:dof_damping_norm}
    \end{figure}

  Comparing \cref{fig:sof_damping_norm} and \cref{fig:dof_damping_norm} we can see that the Luenberger observer control appears to be outperformed by the output feedback control, despite the fact that the former is permitted both current and historical observations and the latter only current observations. This is because the objective of the controls is to minimise the cost \cref{eqn:discrete_cost} rather than flatten the film fastest. Comparing \(t=0\) we can clearly see that the output feedback controls are significantly larger, resulting in a larger cost: \(\sim 11.60\) versus \(\sim 7.39\).

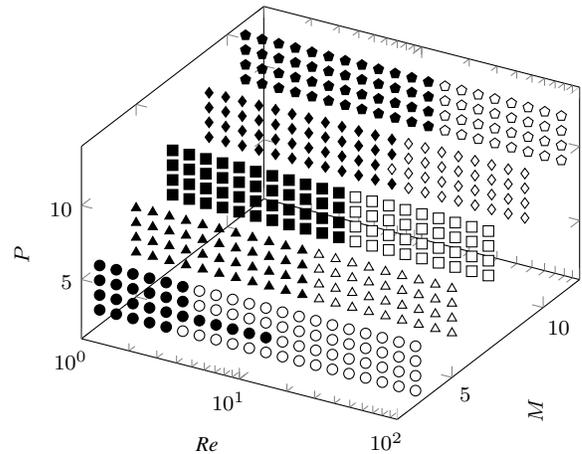
\begin{figure}[tpb]
    \centering
    \begin{tikzpicture}[font=\footnotesize]
      \begin{axis}[
        width=0.95\columnwidth,
        axis lines = box,
        xlabel = {\(\Rey\)},
        ylabel = {\(M\)},
        zlabel = {\(P\)},
        xmin = 1,
        xmax = 100,
        ymin = 2,
        ymax = 12,
        zmin = 1,
        zmax = 14,
        xmode=log,
        view={30}{40},
        ]

        \addplot3 [only marks, fill=black, mark=*,
          y filter/.expression={ y-3 == 0 ? \pgfmathresult : NaN)}]
        table[x=R, y=M, z=P]
        {data_dof-success.dat};

        \addplot3 [only marks, fill=white, mark=*,
          y filter/.expression={ y-3 == 0 ? \pgfmathresult : NaN)}]
        table[x=R, y=M, z=P]
        {data_dof-failure.dat};

        \addplot3 [only marks, fill=black, mark=triangle*,
          y filter/.expression={ y-5 == 0 ? \pgfmathresult : NaN)}]
        table[x=R, y=M, z=P]
        {data_dof-success.dat};

        \addplot3 [only marks, fill=white, mark=triangle*,
          y filter/.expression={ y-5 == 0 ? \pgfmathresult : NaN)}]
        table[x=R, y=M, z=P]
        {data_dof-failure.dat};

        \addplot3 [only marks, fill=black, mark=square*,
          y filter/.expression={ y-7 == 0 ? \pgfmathresult : NaN)}]
        table[x=R, y=M, z=P]
        {data_dof-success.dat};

        \addplot3 [only marks, fill=white, mark=square*,
          y filter/.expression={ y-7 == 0 ? \pgfmathresult : NaN)}]
        table[x=R, y=M, z=P]
        {data_dof-failure.dat};

        \addplot3 [only marks, fill=black, mark=diamond*,
          y filter/.expression={ y-9 == 0 ? \pgfmathresult : NaN)}]
        table[x=R, y=M, z=P]
        {data_dof-success.dat};

        \addplot3 [only marks, fill=white, mark=diamond*,
          y filter/.expression={ y-9 == 0 ? \pgfmathresult : NaN)}]
        table[x=R, y=M, z=P]
        {data_dof-failure.dat};

        \addplot3 [only marks, fill=black, mark=pentagon*,
          y filter/.expression={ y-11 == 0 ? \pgfmathresult : NaN)}]
        table[x=R, y=M, z=P]
        {data_dof-success.dat};

        \addplot3 [only marks, fill=white, mark=pentagon*,
          y filter/.expression={ y-11 == 0 ? \pgfmathresult : NaN)}]
        table[x=R, y=M, z=P]
        {data_dof-failure.dat};
      \end{axis}
    \end{tikzpicture}
    \caption{Successful (in black) and unsuccessful (in white) results for the Luenberger observer method applied to simulations of the Navier-Stokes equations at various \(\Rey\), number of actuators \(M\), and numbers of observers \(P\). The groups of differently shaped points correspond to the values of \(M\): 3, 5, 7, 9, and 11.}%
    \label{fig:dof_damping_scatter}
  \end{figure}

  Although the Luenberger observer control is an obvious improvement over the output feedback control (as can be seen by comparing \cref{fig:sof_damping_scatter,fig:dof_damping_scatter}), unsurprisingly is unable to match the case of full observations at larger Reynolds numbers, although we again see stabilisation above \(\Rey \approx 5\). As for the output feedback LQR problem, as \(\Rey\) grows the nonlinear effects become more important, and so restricting information neglects the energy transferred from large to small modes, which the observer is unable to resolve.

\section{CONCLUSION}
\label{sec:conclusion}
  We demonstrate a pair of techniques to extend a hierarchical control scheme for the Navier-Stokes equations to cases with restricted observations. We find that although such restrictions have a significant impact on the regimes where stabilisation is achievable, the methods are nonetheless able to stabilise a uniform film in a region including and at times exceeding the range of predictive capacity of the underlying simplified models. We note that the Luenberger observer has two major advantages over the output feedback case: it requires significantly simpler and more stable algorithms to compute the controls, which are subsequently able to achieve the same stabilisation results at a lower cost.


\bibliography{references.bib}

\end{document}